\documentclass[letterpaper, 10 pt, conference]{ieeeconf}
\IEEEoverridecommandlockouts
\usepackage{cite}
\usepackage{amsmath,amssymb,amsfonts}
\usepackage{graphicx}

\newtheorem{assumption}{Assumption}

\title{Steady state programming of controlled nonlinear systems via deep dynamic mode decomposition}

\author{Aqib Hasnain, Nibodh Boddupalli, Shara Balakrishnan, and Enoch Yeung
\thanks{A. Hasnain ({\tt\small aqib@ucsb.edu}) and Nibodh Boddupalli are with the Department of Mechanical Engineering, University of California, Santa Barbara \newline \indent S. Balakrishnan is with the Department of Electrical and Computer Engineering, University of California, Santa Barbara \newline \indent E. Yeung ({\tt\small eyeung@ucsb.edu}) is with the Department of Mechanical Engineering, Center for Control, Dynamical Systems, and Computation, and
Biomolecular Science and Engineering, University of California, Santa Barbara}
}

\begin{document}
 
\maketitle

\begin{abstract}
This paper describes the optimal selection of a control policy to program the steady state of controlled nonlinear systems with hyperbolic fixed points. This work is motivated by the field of synthetic biology, in which saddle points are common (along with limit cycles), and the aim is to program cells to perform both digital and analog computation, though developing genetic digital computation has been the main focus. We frame the analog computing challenge of generating a steady state input-output function inside living cells. To program the steady state, a data-driven approach is taken wherein an approximation of the Koopman operator, identified via deep dynamic mode decomposition, is used to describe the dynamics of the system linearly. The new representation of the dynamics are then used to solve an optimization problem for the input which maximizes a direction in state space. Some added structure on the Koopman operator learning process for controlled systems is given for dynamics that are separable in the state and input. Finally, the methods are demonstrated on simulation examples of an incoherent feedforward loop and a combinatorial promoter system, two common network architectures seen in the field of synthetic biology.
\end{abstract}

\section{Introduction}
The field of synthetic biology is concerned with designing and constructing biological modules, biological systems, and biological machines \cite{cameron2014brief} and can be traced back to the 1960s when logic in genetic regulation was discovered \cite{jacob1961genetic}. The key goal is to program living cells to exhibit desired functionality. Biomolecular processes are typically nonlinear, stochastic, and non-modular making the design and construction process difficult. This is where control theory meets synthetic biology to produce robust functionality \cite{del2016control}. 

Thus far, the main focus of the field of synthetic biology has been in designing and building digital logic genetic circuits. A wide variety of digital computation using synthetic genetic circuits has been achieved including but not limited to switches \cite{gardner2000construction}, counters \cite{friedland2009synthetic}, logic gates \cite{singh2014recent}, classifiers \cite{didovyk2015distributed}, edge detectors \cite{tabor2009synthetic}, and most recently digital displays \cite{millacura2019parallel,shin2020programming}. The state-of-the-art software called CELLO \cite{nielsen2016genetic} automates the design of digital logic genetic circuits and has paved the way for larger genetic circuits to be developed. Moreover, network reconstruction algorithms which infer dynamic network architectures have been developed to debug failure modes in these digital logic devices \cite{yeung2015global,ward2009comparison,sontag2004inferring,hasnain2019data}. However, there are cellular mechanisms and limitations that make it unclear whether or not digital computation can scale up to perform complex computations in living cells. 

Biological systems such as colonies of bacteria inherently implement hybrid-analog computing \cite{sarpeshkar2014analog}, which was a staple paradigm of computation in the 1950s and 1960s. Analog computers need to be custom designed for a specific task and are less reliable than their digital counterparts. However, the benefits of analog computing, combined with the non-von Neumann sense of computing, is two-fold, i) it is customizable to address desired functionality (also seen as a disadvantage as noted prior), and ii) it gives rise to a type of computation where memory and processing are collocated, making for a more efficient computing architecture. We want to exploit the analog, tunable nature to perform biological computation. Specifically, we aim to design an input-output function in bacteria that takes chemical inducers as input and outputs the maximum steady state concentration. The bacterial colony as a whole can be viewed as the central processing unit with memory \cite{ben2009learning}. In cellular computing, data may be thought of as being encoded by biomolecules such as DNA strands and molecular biology tools may act on the data to perform various operations \cite{paun2005dna}. 

There exists no general theory for how cells perform computation. Likewise, there is no general theory for many of the \textit{known} biomolecular mechanisms. This provides difficulty when designing and constructing synthetic biological components. Hypothesis-driven modeling is one approach to better understand these biomolecular mechanisms and how to control them \cite{yeung2017biophysical,jayanthi2013retroactivity}. A downside of hypothesis-driven modeling is that the models can be difficult to validate, often requiring iteration upon iteration of experimentation to fit model parameters. Futhermore, when the biological process occurs in new environmental conditions, it is likely that new experimental data will need to be collected for model refinement. 

The above issues and more motivate the need for purely data-driven algorithms to accelerate the advancement of design of genetic circuits. We frame the analog computing challenge - is it possible to harness the natural dynamics of the cell to generate a steady state input-output function? Specifically, how do we utilize time-series measurements of a biological system to design control inputs to achieve a target steady state output from living cells in a data-driven fashion? 

For this, we turn to Koopman operator theory, which is a powerful tool for data-driven analysis of nonlinear dynamical systems. Researchers working in this space have shown that it is possible to identify and learn the fundamental modes of a nonlinear dynamical system from data \cite{rowley2009spectral,mezic2005spectral}. The Koopman operator is an infinite-dimensional linear operator that represents nonlinear dynamics as a dynamically equivalent linear system. The development of Dynamic Mode Decomposition (DMD) \cite{schmid2010dynamic} has led to rapid growth in the use of Koopman spectral analysis of nonlinear dynamical systems in areas such as system identification \cite{mauroy2019koopman,boddupalli2019koopman}, prediction and control \cite{korda2018linear,korda2019optimal,proctor2016dynamic,proctor2018generalizing}, and sensor placement \cite{sinha2016operator,hasnain2019optimal}. More recently, learning higher dimensional Koopman operators from data has become computationally tractable, largely due to advances in integrating machine learning and deep learning to generate efficient representations of observable bases \cite{yeung2019learning,lusch2018deep,otto2019linearly,takeishi2017learning}.

Previously, researchers have explored control strategies to reprogram the steady states of cooperative monotone
dynamical systems \cite{shah2018reprogramming} and in \cite{del2017blueprint}, a synthetic genetic
feedback controller that dynamically steers the concentration
of a genetic regulatory network’s key transcription factors
was developed to reprogram the steady state. We propose to
alter the steady state, or generate a genetic program with output being steady state values, of nonlinear systems with hyperbolic
fixed points by designing an optimal control policy for the
Koopman model identified directly from data. Specifically,
we utilize deep Dynamic Mode Decomposition (deepDMD)
[33] to learn approximate Koopman invariant subspaces
of dynamical systems to program their steady state. The
framework can be extended to systems with other types of
attractors and can be modified to solve many optimal control
problems e.g. target state and reference tracking.

In section \ref{sec:koop}, we briefly review Koopman operator theory for discrete-time nonlinear systems the numerical approximations of the Koopman operator. Section \ref{sec:ssProgramming} describes the steady state programming framework and discusses issues
that arise when dealing with mixed terms of the state and
input. Lastly, we demonstrate our method on two nonlinear
systems that are commonly seen in biomolecular feedback
systems and synthetic biology.

\section{Koopman Operators and their finite-dimensional approximation} \label{sec:koop}
In this section we briefly review Koopman operator theory as it applies to controlled systems. See \cite{kaiser2019data} for an overview of the field. In this paper, we consider discrete-time dynamical systems to be consistent with the nature of time-series data. 

\subsection{Koopman operator theory for control systems}
Consider a non-affine control discrete-time nonlinear dynamical system of the form
\begin{equation}
    x_{k+1} = f(x_k,u_k),
    \label{eq:stateDyn}
\end{equation}
where $x_k \in \mathbb{R}^n$ is the state of the system, $u_k \in \mathbb{R}^m$ is the control input, and $f: \mathbb{R}^n \oplus \mathbb{R}^m \rightarrow \mathbb{R}^n$ is the analytic and unknown transition mapping. The dynamics can be always be decomposed in the following way
\begin{equation}
    x_{k+1} = f_x(x_k) + f_{xu}(x_k,u_k) + f_u(u_k),
    \label{eq:stateDynDecomp}
\end{equation}
where $f_x$, $f_{xu}$, $f_u$ are terms only in $x$, mixed terms in $x$ and $u$, and only in $u$, respectively. The Koopman operator acts on a set of observables and we consider each observable as an element of an infinite-dimensional Hilbert space $\mathcal{F}$. 

The observables $\psi: \mathbb{R}^n \oplus \mathbb{R}^m \rightarrow \mathbb{R}$ are functions of the state and the input such that they can be functions of only $x$, mixed terms of both $x$ and $u$, and only of $u$. The observables can also be thought of as vector-valued functions of the state and the input such that $\psi: \mathbb{R}^n \oplus \mathbb{R}^m \rightarrow \mathbb{R}^{o}$.

Then the Koopman operator, $\mathcal{K}: \mathcal{F} \rightarrow \mathcal{F}$ acts on the Hilbert space of observables to produce linear dynamics
\begin{equation}
    \mathcal{K}\psi(x_k,u_k) \triangleq \psi(f(x_k,u_k),u_{k+1}) 
    \label{eq:koopEqgeneral}
\end{equation}
so that
\begin{equation}
    \psi(x_{k+1},u_{k+1}) = \mathcal{K}\psi(x_k,u_k)
\end{equation}
which can be decomposed similarly to (\ref{eq:stateDynDecomp}) into the three components
\begin{equation}
    \begin{bmatrix}\psi_x(x_{k+1}) \\ \psi_{xu}(x_{k+1},u_{k+1}) \\ \psi_u(u_{k+1}) \end{bmatrix} = \begin{bmatrix} \mathcal{K}_x & \mathcal{K}_{xu} & \mathcal{K}_u \\\mathcal{K}_{21} & \mathcal{K}_{22} & \mathcal{K}_{23} \\ \mathcal{K}_{31} & \mathcal{K}_{32} & \mathcal{K}_{33} \end{bmatrix} \begin{bmatrix}\psi_x(x_{k}) \\ \psi_{xu}(x_{k},u_{k}) \\ \psi_u(u_{k}) \end{bmatrix}
    \label{eq:koopDecomp}
\end{equation}
where $\psi_x: \mathbb{R}^n \rightarrow \mathbb{R}^{n_L}$, $\psi_{xu}: \mathbb{R}^n \oplus \mathbb{R}^m \rightarrow \mathbb{R}^{M_L}$, and $\psi_x: \mathbb{R}^m \rightarrow \mathbb{R}^{m_L}$. 
\begin{assumption}
\textit{We assume that the inputs $u_k$ are exogenous disturbances without state space dynamics. }
\label{ass:exogenousInputs}
\end{assumption}
The above Assumption \ref{ass:exogenousInputs} allows us to be concerned only with the first block equation in (\ref{eq:koopDecomp}), i.e.
\begin{equation}
    \psi_x(x_{k+1}) = \mathcal{K}_x \psi_x(x_{k}) + \mathcal{K}_{xu}\psi_{xu}(x_{k},u_{k}) + \mathcal{K}_u\psi_u(u_{k}). 
    \label{eq:xkoopDecomp}
\end{equation}
We refer the readers to \cite{yeung2018koopman} for more discussion on this decomposition and a proof for the existence of Koopman operators for nonlinear systems. 
\subsection{Finite dimensional approximations}
\subsubsection{Dynamic mode decomposition with control}
Dynamic mode decomposition with control (DMDc) developed in \cite{proctor2016dynamic,proctor2018generalizing} is an extension of DMD for systems with external actuation and is used to identify a finite dimensional approximation to the Koopman operator, $K$. The idea is to find the best-fit linear operators $A$ and $B$ to provide the following relationship:
\begin{equation*}
    x_{k+1} = Ax_k + Bu_k
\end{equation*}
for measurement $x_k$, present control $u_k$, and future measurement $x_{k+1}$. To do this, data snapshots are collected from the dynamical system (\ref{eq:stateDyn}) of the state and input over time and formed into the following data matrices:
\begin{equation*}
\begin{aligned}
    X_p &= \begin{bmatrix} x_0 & x_1 & \hdots & x_{N-1} \end{bmatrix},\\
    X_f &= \begin{bmatrix} x_1 & x_2 & \hdots & x_N \end{bmatrix},\\
    U_p &= \begin{bmatrix} u_0 & u_1 & \hdots & u_{N-1} \end{bmatrix},
\end{aligned}
\end{equation*}
where N is the number of snapshots collected. The above matrices are known as snapshot matrices since each column represents a snapshot at time $k$ of the state of the dynamical system. If data is collected from more than one trajectory, the additional time-series can be appended as additional columns in the snapshot matrices. To identify the best-fit linear operators, the algorithm solves the following optimization problem:
\begin{equation*}
\begin{aligned}
    \begin{bmatrix}A&B\end{bmatrix} &= \text{arg}\min_{A,B} \sum_{i=0}^{N-1} || x_{i+1} - Ax_i-Bu_i  ||_2,  \\
    &= \text{arg}\min_{A,B} || X_f - AX_p - BU_p||_F, \\
    & = X_f\begin{bmatrix}X_p \\ U_p\end{bmatrix}^{\dagger}
\end{aligned}
\end{equation*}
where $\dagger$ represents the Moore-Penrose pseudinverse. Here $A$ corresponds to $K_x$ and $B$ to $K_u$ while $K_{xu}=0$. Thus DMDc is limited in its system identification capabilites as it does not account for mixed terms of $x$ and $u$ in its formulation. 
If trajectories of the dynamical system without inputs is available, the state transition matrix $A$ can first be identified and treated to be a constant in the learning of the control matrix $B$ from the trajectories with inputs \cite{hasnain2019data}. This disambiguates the impact of the control input on the state dynamics. 

\subsubsection{Deep dynamic mode decomposition with control}
Another approach to learning Koopman invariant subspaces for controlled nonlinear systems is given in \cite{williams2016extending}, but a difficulty in this approach is that the user must choose the observables manually. We use the deep learning approach deep Dynamic Mode Decomposition (deepDMD) developed in \cite{yeung2019learning} to identify an efficient representation of the invariant subspace. The approach here is to represent the observable functions as many compositions of linear transformations and ReLU (Rectified Linear Units) functions and learn both the observables and the Koopman operator simultaneously. The neural network is tasked with solving the following optimization problem
\begin{equation*}
\begin{aligned}
    &\min_{K_x,K_{xu},K_u,\theta} || \Psi_x(X_f,\theta) - K_x\Psi_x(X_p,\theta) \\& \qquad \qquad \quad - K_{xu}\Psi_{xu}(X_p,U_p,\theta) - K_u\Psi_{u}(U_p,\theta)  ||_F \\& \qquad \qquad \quad
    + \lambda_1||K||_2 + \lambda_2||\theta||_1
\end{aligned}
\end{equation*}
where the parameters $\theta$ are the neural networks biases and weights and $\lambda_1$ and $\lambda_2$ are regularization parameters. The snapshot matrix $\Psi_x$, is now represented as
\begin{equation*}
    \Psi_x(X_p,\theta) = \begin{bmatrix} \psi_x(x_0,\theta) & \psi_x(x_1,\theta) & \hdots & \psi_x(x_{N-1},\theta) \end{bmatrix}
\end{equation*} 
and the other snapshot matrices are similarly formed. The neural network identifies a model that relates the data in observable space as
\begin{equation*}
    \psi_x(x_{k+1}) = K_x\psi_x(x_k) + K_{xu}\psi_{xu}(x,u) + K_u\psi_u(u),
\end{equation*}
where we have dropped the dependency on $\theta$ for brevity.  

\section{Steady state programming} \label{sec:ssProgramming}

In this section, when we refer to the Koopman operator, we are referring to the finite-dimensional approximation of the actual Koopman operator. The focus of this paper is to design an input $u$ to maximize a single direction of the state at equilibrium. We emphasize that we are not solving a dynamic control policy planning problem, but that the control policy is a constant input that maximizes the equilibrium value of a single direction in state space. Furthermore, we consider nonlinear systems which have hyperbolic fixed points. At steady state we have the optimization problem 
\begin{equation}
\begin{aligned}
&\min_{u \in \mathcal{U}} \quad - \hat{e}^{\top}_i x_{k,e} \\
& \quad \textrm{s.t.} \quad  x_{k,e} = f(x_{k,e},u_k), 
  \label{eq:ssOpt}
\end{aligned}
\end{equation}
where $\mathcal{U}$ is a bounded set of inputs, $\hat{e}_i$ is the unit column vector in the $i^{th}$ direction of the state (the state we aim to maximize), and the subscript $e$ denotes equilibrium values. Since the form of $f$ is unknown, it is difficult to design an input which solves this program. Therefore, we propose to identify a Koopman invariant subspace by extending the states and the inputs into a space of observables. The identified Koopman model then allows us to formulate the optimization problem in the lifted space where the problem has some tractability.

\subsection{Steady state programming using DMD with linear control}

For a linear time-invariant system with $(I - A)$ invertible, if $u$ does not come from a bounded set $\mathcal{U}$ but from an unbounded set $\hat{\mathcal{U}}$ and the observable mapping is the identity, there is no bounded $u$ that solves optimization problem (\ref{eq:ssOpt}). To see this, under the conditions above, the optimization problem (\ref{eq:ssOpt}) becomes
\begin{equation*}
  \min_{u \in \hat{\mathcal{U}}} \quad -\hat{e}^{\top}_i(I - A)^{-1}Bu_k
\end{equation*}
and note that the objective is linear in $u$. Differentiating with respect to $u$ we have
\begin{equation*}
    \frac{\partial}{\partial u_k} -\hat{e}^{\top}_i(I - A)^{-1}B u_k = -\hat{e}^{\top}_i(I - A)^{-1}B
\end{equation*}
and so for any $(A,B)$, there is no bounded $u^*$ that solves the optimization problem. If $u$ comes from the bounded set $\mathcal{U}$, then the minimum value of the objective on $\mathcal{U}$ occurs on the boundary $\partial \mathcal{U}$ due to the maximum modulus principle. This results from the fact that the objective function to be optimized is linear in $u$ and therefore it is an uninteresting case. This case corresponds to identifying a linear model using DMD with control, where the observable mapping is simply the identity and the measurement snapshots are all that is used to compute the state-transition matrix and the control matrix. For nonlinear systems that can be represented linearly by considering a linear basis, the above approach will work effectively. However, we aim to develop a method which can handle nonlinear systems that \textit{require} a basis that is a nonlinear function of the original states. 

\subsection{Steady state programming using deepDMD with nonlinear control }
Consider the case where the observable mapping is not the identity.
\begin{assumption}
\textit{We assume that the state transition Koopman operator $K_x$ has $n_L$ linearly independent eigenvectors and that none of the corresponding $n_L$ eigenvalues are equal to 1. }
\label{ass:invertibleIminusKx}
\end{assumption}
At equilibrium, (\ref{eq:xkoopDecomp}) becomes
\begin{equation}
    \psi_x(x_{e}) = K_x\psi_x(x_{e}) + K_{xu}\psi_{xu}(x_{e},u_{k}) + K_u\psi_u(u_{k})
    \label{eq:ssKoopEq}
\end{equation}
and under Assumption \ref{ass:invertibleIminusKx}, $(I - K_x)$ is invertible, so we have
\begin{equation*}
    \psi_x(x_{e}) = (I - K_x)^{-1} \left[K_{xu}\psi_{xu}(x_{e},u_{k}) + K_u\psi_u(u_{k})\right].
\end{equation*}
For the optimization problem (\ref{eq:ssOpt}) we now have
\begin{equation}
\begin{aligned}
&\min_{\psi_u(u_k) \in \mathcal{U}_{\Psi}} \quad - \hat{e}^{\top}_i \psi_x(x_{e}) \\
& \quad \textrm{s.t.} \quad \psi_x(x_{e}) = (I - K_x)^{-1} \\ & \qquad \qquad \left[K_{xu}\psi_{xu}(x_{e},u_{k}) + K_u\psi_u(u_{k})\right] 
\end{aligned}
\label{eq:ssOptKoop}
\end{equation}
stating that we want to find the input $\psi_u(u)$ that maximizes the equilibrium value of a single direction in observable space. Here $\mathcal{U}_{\Psi}$ is a bounded set of lifted inputs. If we assume away the mixed terms of $x$ and $u$ then the problem becomes immediately tractable as we have $\psi_x(x_e)$ as a function of just terms in $\psi_u(u)$. However, assuming away the mixed terms seemingly requires a strong assumption about the nonlinearities present in original dynamics. Specifically it assumes that the input affects the state dynamics through additive means \textit{only}. This is typically not the case in biological systems. Another important point to make is that this is an implicit optimization of $\psi_x(x_e)$, another layer of difficulty when dealing with the mixed terms.

Since we transform the state space coordinates to an abstract observable-space, we want to keep the interpretability of the original system intact. Therefore we make the following assumption:
\begin{assumption}
\textit{We assume state-inclusive and input-inclusive observables $\psi_{x}(x)$ and $\psi_u(u)$, respectively, i.e.}
\begin{equation}
\begin{aligned}
    \psi_x(x) &= \begin{bmatrix} x \\ \varphi_x(x) \end{bmatrix}, \\ 
    \psi_u(u) & = \begin{bmatrix} u \\ \varphi_u(u) \end{bmatrix}
\end{aligned}
\end{equation}
\textit{with $\varphi_x(x) \in \mathbb{R}^{n_L-n}$ and $\varphi_u(u) \in \mathbb{R}^{m_L -m}$.}
\end{assumption}
Again, since the mixed terms of $x$ and $u$ cause issues in solving program (\ref{eq:ssOptKoop}), we want to examine the structure of these terms when an EDMD \cite{williams2016extending} with control approach  is used. A typical EDMD dictionary may consist of functions such as monomials, polynomials (Legendre and Hermite), radial basis functions, trigonometric functions, and logistic functions. It can be shown that for most of these common dictionary elements, the mixed terms are separable in $x$ and $u$. To illustrate what we want to achieve, consider an example system with state-inclusive ($n=2$ and $m=2$) observables with up to second order monomials in the state and input. Then for $\psi_x: \mathbb{R}^n \rightarrow \mathbb{R}^{n_L = 5}$ and $\psi_u: \mathbb{R}^m \rightarrow \mathbb{R}^{m_L = 5}$ we have that
\begin{equation*}
    \psi_x(x) = \begin{bmatrix} x_1 \\ x_2 \\ x_1x_2 \\ x_1^2 \\ x_2^2 \end{bmatrix}, \quad
    \psi_u(u) = \begin{bmatrix} u_1 \\ u_2 \\ u_1u_2 \\ u_1^2 \\ u_2^2 \end{bmatrix}
\end{equation*}
where we have left out the zeroth order terms for brevity. For this example, the subscripts on the states and inputs represent different components of the state and input, not timepoints. As we have no knowledge of the form of the dynamics $f$, we may consider forming all possible combinations of $\psi_x(x)$ and $\psi_u(u)$ and using those elements to form $\psi_{xu}(x,u)$, 
\begin{equation*}
    \psi_{xu}(x,u) = \begin{bmatrix} x_1\psi_u(u) \\ x_2\psi_u(u) \\ x_1x_2\psi_u(u) \\ x_1^2\psi_u(u) \\ x_2^2\psi_u(u) \end{bmatrix} \text{or} \quad \psi_{xu}(x,u) = \begin{bmatrix} u_1\psi_x(x) \\ u_2\psi_x(x) \\ u_1u_2\psi_x(x) \\ u_1^2\psi_x(x) \\ u_2^2\psi_x(x) \end{bmatrix}
\end{equation*}
with $\psi_{xu}: \mathbb{R}^n \oplus \mathbb{R}^m \rightarrow \mathbb{R}^{M_L = 25}$, such that $M_L = n_L \times m_L$. Note that $\psi_{xu}(x,u)$ can now be written as a matrix times $\psi_u(u)$ or a matrix times $\psi_x(x)$ as
\begin{equation}
    \psi_{xu}(x,u) = \begin{bmatrix} D_x(x_1) \\ D_x(x_2) \\ D_x(x_1x_2) \\ D_x(x_1^2) \\ D_x(x_2^2) \end{bmatrix} \psi_u(u) = M_x \psi_u(u)
    \label{eq:sepxu}
\end{equation}
or 
\begin{equation}
    \psi_{xu}(x,u) = \begin{bmatrix} D_u(u_1) \\ D_u(u_2) \\ D_u(u_1u_2) \\ D_u(u_1^2) \\ D_u(u_2^2) \end{bmatrix} \psi_x(x) = M_u \psi_x(x)
    \label{eq:sepux}
\end{equation}
where $D_x(\cdot) \in \mathbb{R}^{n_L\times m_L}$ and $D_u(\cdot) \in \mathbb{R}^{m_L\times n_L}$ are diagonal matrices with their argument as the diagonal elements. The matrices $M_x \in \mathbb{R}^{M_L \times m_L}$ and $M_u \in \mathbb{R}^{M_L \times n_L}$ consist of terms only in $x$ and only in $u$, respectively.

A more general statement is that if $f(x,u)$ in (\ref{eq:stateDynDecomp}) is separable into some $G(u)h(x)$, then the mixed observable $\psi_{xu}(x,u)$ are also similarly separable. To see this, consider continuous-time dynamics $\dot{x} = f(x,u)$ for ease of presentation. If $f(x,u)$ is multiplicatively separable in $x$ and $u$, it can be written as 
\begin{equation}
    f(x,u) = WG(u)h(x),
\end{equation}
with the weight matrix $W \in \mathbb{R}^{n\times M_L}$, the matrix $G(u) \in \mathbb{R}^{M_L \times n_L}$ is a function of $u$, and the vector $h(x) \in \mathbb{R}^{n_L}$. From the Koopman generator \cite{budivsic2012applied}, we have that 
\begin{equation}
    \dot{\psi}_x(x) = \frac{\partial \psi_x(x)}{\partial x }f(x,u) = \frac{\partial \psi_x(x)}{\partial x } WG(u)h(x).
\end{equation}
and therefore, any mixed terms are entirely separable in the observable space as well. This may mean that the deep learning approach can take care of the necessity of having mixed terms without explicitly learning the $K_{xu}$ parameters. 

Using the idea that we can separate the mixed terms as in (\ref{eq:sepxu}) and (\ref{eq:sepux}), (\ref{eq:ssKoopEq}) can be written as
\begin{equation}
    \psi_x(x_{e}) = K_x\psi_x(x_{e}) + K_{xu}M_u \psi_x(x_e) + K_u\psi_u(u_{k})
    \label{eq:ssKoopEqsepux}
\end{equation}
or as 
\begin{equation}
    \psi_x(x_{e}) = K_x\psi_x(x_{e}) + K_{xu}M_x \psi_u(u_k) + K_u\psi_u(u_{k}).
    \label{eq:ssKoopEqsepxu}
\end{equation}
From (\ref{eq:ssKoopEqsepux}) we have that at equilibrium, 
\begin{equation*}
    \psi_x(x_{e}) = (I-K_x-K_{xu}M_u)^{-1}K_u\psi_u(u_{k})
\end{equation*}
if $(I-K_x-K_{xu}M_u)$ is invertible. From (\ref{eq:ssKoopEqsepxu}) we have that at equilibrium,
\begin{equation*}
    \psi_x(x_{e}) = (I-K_x)^{-1}\left[K_{xu}M_x+K_u\psi_u(u_{k})\right]
\end{equation*}
and from Assumption \ref{ass:invertibleIminusKx} $(I-K_x)$ is invertible. 

Of these two options, if we use (\ref{eq:ssKoopEqsepxu}) to program the steady state of nonlinear systems, it only requires that $(I-K_x)$ to be invertible, which has already been assumed. Equation (\ref{eq:ssKoopEqsepux}) on the other hand requires that $(I-K_x-K_{xu}M_u)$ be invertible. An advantage is that the right hand side is entirely a function of $u$ whereas (\ref{eq:ssKoopEqsepxu}) has a right hand side which is a function of both $u$ and $x$, requiring implicit optimization. From this analysis, we can again rewrite the steady state programming problem (\ref{eq:ssOptKoop}) as 
\begin{equation}
\begin{aligned}
&\min_{\psi_u(u_k) \in \mathcal{U}_{\Psi}} \quad - \hat{e}^{\top}_i \psi_x(x_{e}) \\
& \quad \textrm{s.t.} \quad \psi_x(x_{e}) = (I-K_x-K_{xu}M_u)^{-1}K_u\psi_u(u_{k})  \\
& \quad \textrm{or s.t.} \quad  \psi_x(x_{e}) = (I-K_x)^{-1}\left[K_{xu}M_x+K_u\psi_u(u_{k})\right]  \\
& \quad \textrm{or s.t.} \quad \psi_x(x_{e}) = (I-K_x)^{-1}K_u\psi_u(u_{k})
\label{eq:ssOptFinal}
\end{aligned}
\end{equation}
where the last constraint comes from setting $K_{xu}=0$ i.e. stating that $\psi_{xu}$ terms are unimportant. We consider this last case since deepDMD offers the flexibility to learn accurate Koopman invariant subspaces even without considering the mixed terms, as will be demonstrated. 

This program is in its most general form nonlinear and will yield local optima. EDMD or any neural network approach do not provide a structure to guarantee a global optimum for this problem and is thus NP-hard. We use the Sequential Least Squares Quadratic Programming (SLSQP) solver from the scipy.optimize python package to solve program (\ref{eq:ssOptFinal}). If DMD is used for steady state programming, any convex optimization solver will suffice. 

\section{Numerical Results}
\subsection{Example: Feedforward Loop}
We consider the following nonlinear dynamical system describing a feedforward loop of five proteins (states) under the influence of two inducers (inputs):

\begin{equation} \label{eq:iffl}
\begin{aligned}
    \dot{x}_0  &= \frac{k_0u_0}{1+\frac{u_1}{K_{d_4}}} - \delta_0x_0 \\
    \dot{x}_1  &= \frac{k_1u_1}{1+\frac{x_0}{K_{d_1}}} - \delta_1x_1 \\
    \dot{x}_2  &= k_2x_1+k_3u_0 - \delta_2x_2 \\
    \dot{x}_3  &= \frac{k_4u_1}{1+\frac{x_2}{K_{d_2}}} - \delta_3x_3 \\
    \dot{x}_4 &= \frac{k_5u_0}{1+\frac{x_3}{K_{d_3}}} - \delta_4x_4 .
\end{aligned}
\end{equation}

\begin{figure}
\centering
\includegraphics[scale=0.3]{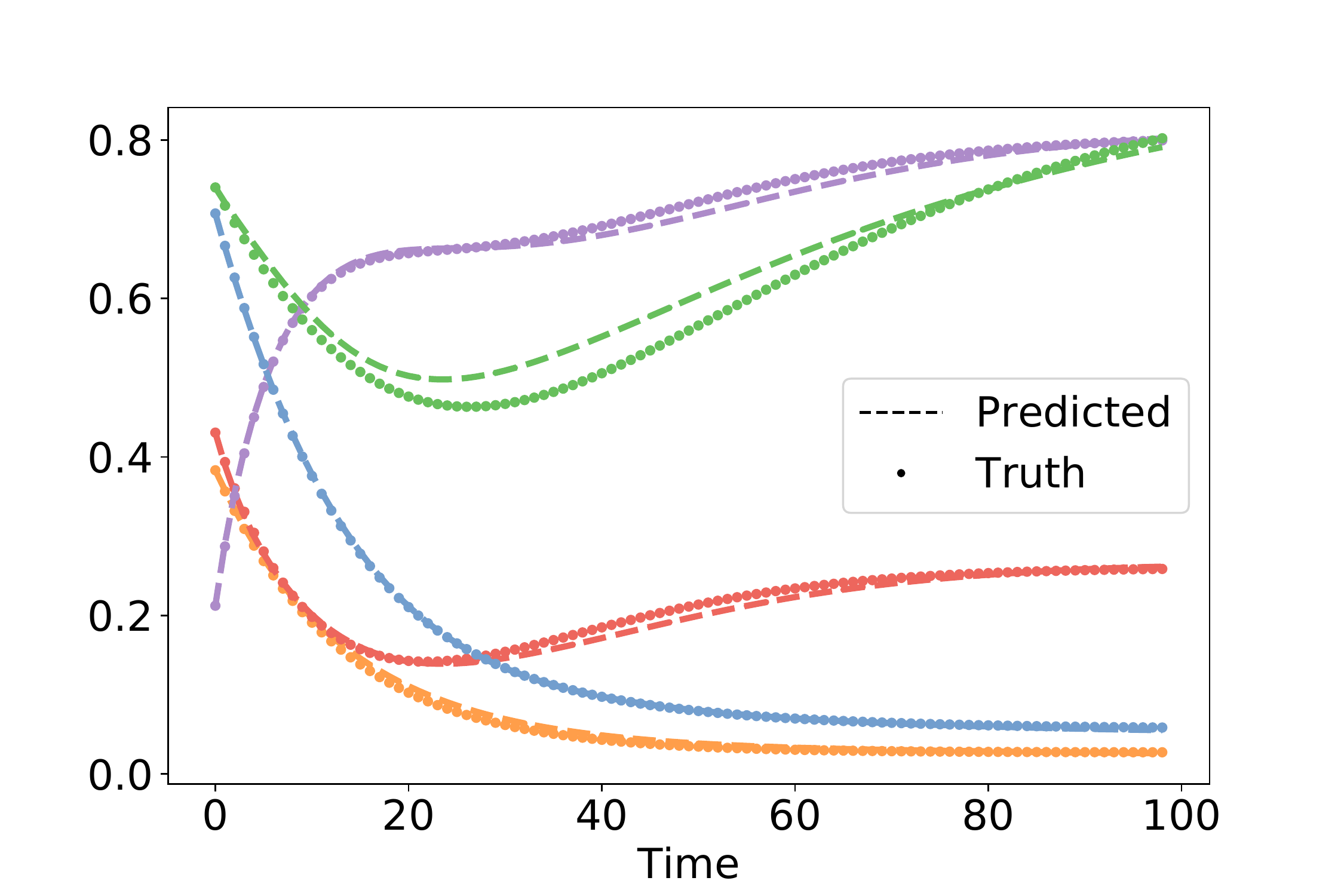}
\caption{99 step prediction from a new initial condition using the neural network identified Koopman model of the feedforward loop (\ref{eq:iffl}). Predictions are given by dashed lines while the true trajectory is given by dotted lines.}
\label{fig:ifflPred}
\end{figure}

\begin{figure}
\centering
\hspace*{-0.3cm} 
\includegraphics[scale=0.25]{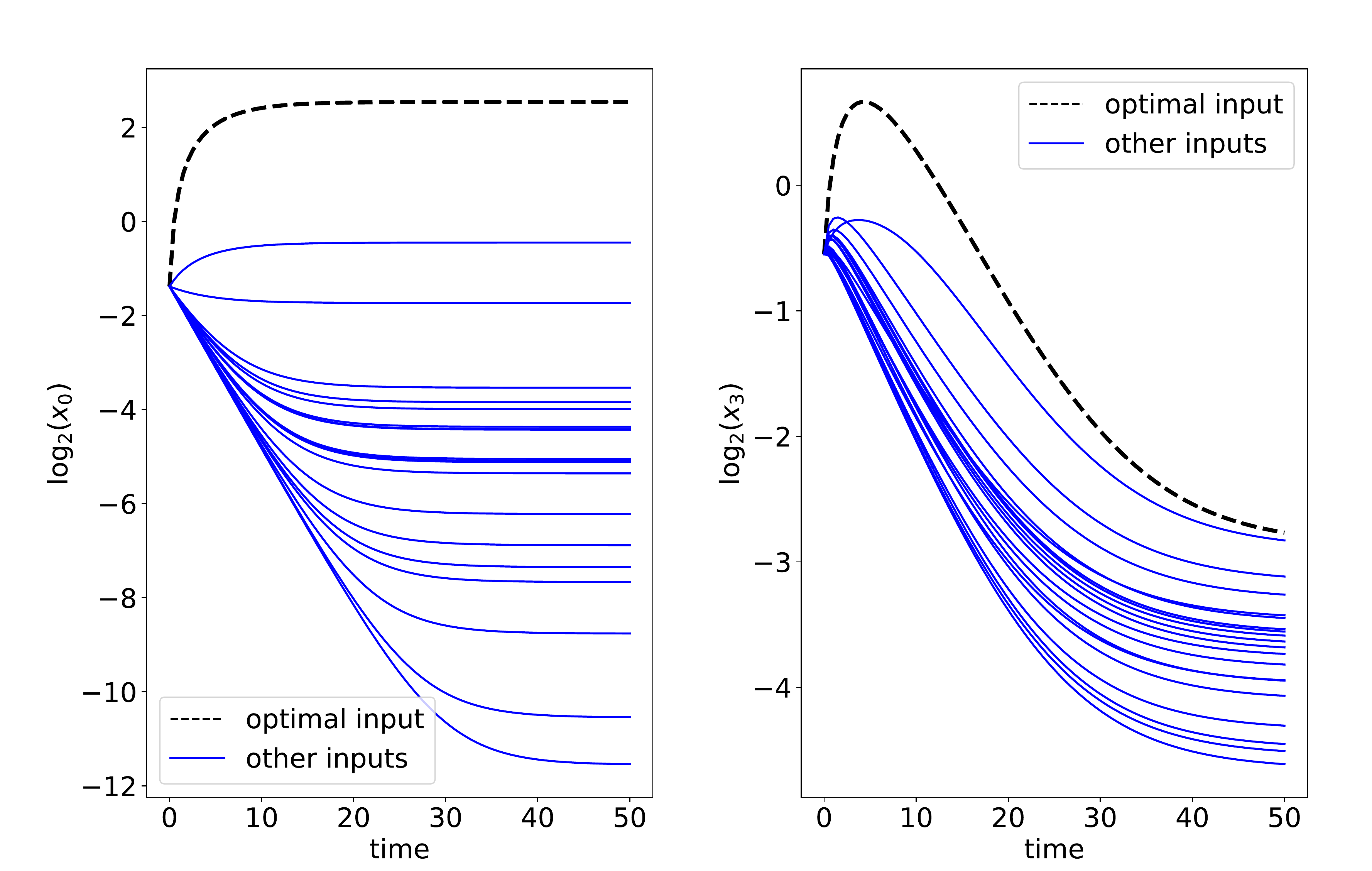}
\caption{Results of the optimal input computed from the steady state programming problem applied to states $x_0$ (left) and $x_3$ (right) of the feedforward loop (\ref{eq:iffl}). The trajectory corresponding to the optimal input is given in black dashed lines and the other trajectories correspond to other inputs that are suboptimal.}
\label{fig:iffl_optInput}
\end{figure}

\normalsize
The states, $x$, describe the concentration of proteins, the inputs $u$ are inducers which activate or repress the protein production and $K_d, \delta, k$ are constant parameters. For simulation, a fourth-order Runge-Kutta scheme is used. 100 timesteps from 100 different initial conditions each with a different input are used for training and testing. The inputs are always step inputs in these simulations. The neural network is a dense feedforward network with the output being the observable basis which spans the approximate Koopman invariant subspace. We take the case where the mixed terms are assumed to be zero. To evaluate accuracy of our Koopman model, we perform a multi-step prediction test. Predictions on a single test trajectory (a trajectory that the neural network has not seen) is given in Figure \ref{fig:ifflPred}. This prediction is done starting from an initial condition and predicting 99 steps into the future and compared with the actual trajectory. The error for this trajectory is $\approx 5\%$ and similar errors are computed for 49 other test trajectories. We then solve program (\ref{eq:ssOptFinal}) and obtain the optimal step input ($0 < u^* < 10$) to achieve a maximum steady state for both state $x_0$ and state $x_3$ as can be seen in Figure \ref{fig:iffl_optInput}. Each trajectory starts from the same initial condition, but only the trajectory corresponding to the black dashed line sees the optimal input and all the other trajectories see a random input between ($0 < u < 10$). It can be seen that the optimal input computed from our framework does give the maximum steady state when applied to the system. Note that the optimal input was applied to the original nonlinear system, \textit{not} to the identified Koopman model.

\subsection{Example: Combinatorial promoter}
For a second example, we consider a combinatorial promoter which takes both a repressor and an activator as input, allowing genes to be switched on and off based on the concentration of the activator and repressor \cite{del2015biomolecular}. The dynamics are given by the following set of differential equations

\begin{equation} \label{eq:combPromoter}
\begin{aligned}
    \dot{x}_0  &= -k_{1f}x_0u_0 + k_{1r}x_2 \\
    \dot{x}_1  &= -k_{2f}x_1u_1 + k_{2r}x_3 -k_{4f}x_1x_4 + k_{4r}x_6 - k_{5f}x_1x_5 \\ &+ k_{5r}x_7 + 0.2x_{10} \\
    \dot{x}_2  &= k_{1f}x_0u_0 - k_{1r}x_2 - k_{3f}x_2x_4 + k_{3r}x_5 - k_{6f}x_2x_6 \\ &+ k_{6r}x_7\\
    \dot{x}_3  &= k_{2f}x_1u_1 - k_{2r}x_3 - 0.2x_3 \\
    \dot{x}_4  &= -k_{3f}x_2x_4 + k_{3r}x_5 - k_{4f}x_1x_4 + k_{4r}x_6\\
    \dot{x}_5  &= k_{3f}x_2x_4 - k_{3r}x_5 - k_{5f}x_1x_5 + k_{5r}x_7 - k_{7f}x_5x_8 \\ &+ k_{7r}x_9 + k_{8f}x_9\\
    \dot{x}_6  &=  -k_{6f}x_2x_6 + k_{6r}x_7 - k_{4r}x_6 + k_{4f}x_1x_4\\
    \dot{x}_7  &=  k_{5f}x_1x_5 - k_{5r}x_7 + k_{6f}x_2x_6 - k_{6r}x_7\\
    \dot{x}_8  &=  -k_{7f}x_5x_8 + (k_{7r}+ k_{8f})x_9\\
    \dot{x}_9  &=  k_{7f}x_5x_8 - (k_{7r}+ k_{8f})x_9\\
    \dot{x}_{10}  &=  k_{8f}x_9 - \delta x_9^2\\
\end{aligned}
\end{equation}

\begin{figure}
\centering
\hspace*{-0.4cm} 
\includegraphics[scale=0.26]{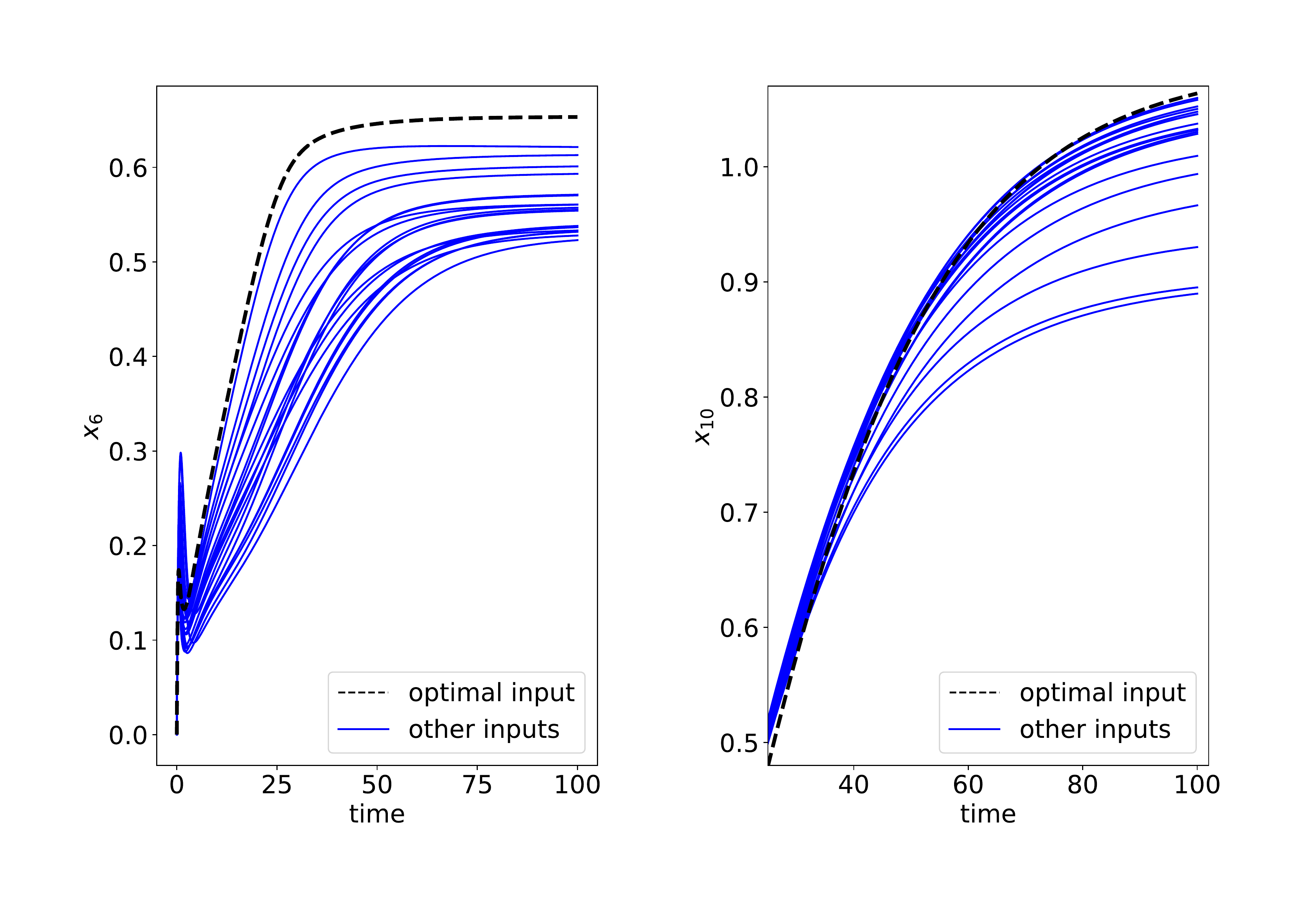}
\caption{Results of the optimal input computed from the steady state programming problem applied to states $x_6$ (left) and $x_{10}$ (right) of the combinatorial promoter system (\ref{eq:combPromoter}). The trajectory corresponding to the optimal input is given in black dashed lines and the other trajectories correspond to other inputs that are suboptimal.}
\label{fig:combPromoter_optInput}
\end{figure}

\normalsize
where the inputs $u$ are the repressor or activator and are continuous, monotonically increasing functions of time. For simulation, a fourth-order Runge-Kutta scheme is used. 1000 timesteps from 80 different initial conditions each with a different input are used for training and testing. The neural network is the same as in Example 1. We again take the case where the mixed terms are assumed to be zero since the deep learning approach is flexible enough to provide a model that predicts well many steps into the future. We again solve program (\ref{eq:ssOptFinal}) and obtain the optimal step input ($0.01 < u^* < 1$) to achieve a maximum steady state for both state $x_6$ and state $x_{10}$ as can be seen in Figure \ref{fig:combPromoter_optInput}. Each trajectory starts from the same initial condition, but only the trajectory corresponding to the black dashed line sees the optimal input and all the other trajectories see a random input between ($0.01 < u < 1$). For this example of a combinatorial promoter, the learned optimal input gives the maximum steady state value for each direction in state space. Again, we want to emphasize that the optimal input was identified using the Koopman model, but it was applied to the original nonlinear system.  


These numerical results demonstrate that it is possible, by harnessing the natural analog computational power of cells and through a data-driven operator theoretic framework, to generate steady state input-output functions. 

\section{Conclusions}
In this paper we presented a formulation for programming the steady state of controlled nonlinear systems with hyperbolic fixed points. We used deep dynamic mode decomposition (deepDMD) with control to compute approximate Koopman invariant subspaces and Koopman operators which represent the original nonlinear system as an approximate linear system. This allowed us to pose and solve an optimization problem wherein we maximize the steady state value of a single direction in state space. The formulation can be extended to handle various cases e.g. maximize the ratio of two directions in state space at steady state. The formulation can be slightly revised to solve common optimal control tasks e.g. reaching a target state and dynamic reference tracking. The method can also be extended to nonlinear systems with other types of attractors. It was shown that mixed terms of the state and the input can provide difficulties and we discuss a way of dealing with these terms. Finally, we demonstrated our method on two example nonlinear systems that are commonly dealt with in biological processes and briefly discuss broader biological computation implications.  

\section*{Acknowledgements}
The authors thank Igor Mezic, Robert Egbert, Bassam Bamieh, Sai Pushpak, Sean Warnick, and Umesh Vaidya for stimulating conversations. This work was supported by a Defense Advanced Research Projects Agency (DARPA) Grant No. DEAC0576RL01830 and an Institute of Collaborative Biotechnologies Grant. Any opinions, findings and conclusions or recommendations expressed in this material are those of the author(s) and do not necessarily reflect the views of the Defense Advanced Research Projects Agency (DARPA), the Department of Defense, or the United States Government.

\bibliographystyle{abbrv}
\bibliography{main}

\begin{thebibliography}{10}

\bibitem{ben2009learning}
E.~Ben-Jacob.
\newblock Learning from bacteria about natural information processing.
\newblock {\em Annals of the New York Academy of Sciences}, 1178(1):78--90,
  2009.

\bibitem{boddupalli2019koopman}
N.~Boddupalli, A.~Hasnain, S.~P. Nandanoori, and E.~Yeung.
\newblock Koopman operators for generalized persistence of excitation
  conditions for nonlinear systems.
\newblock In {\em 2019 IEEE 58th Conference on Decision and Control (CDC)},
  pages 8106--8111. IEEE, 2019.

\bibitem{budivsic2012applied}
M.~Budi{\v{s}}i{\'c}, R.~Mohr, and I.~Mezi{\'c}.
\newblock Applied koopmanism.
\newblock {\em Chaos: An Interdisciplinary Journal of Nonlinear Science},
  22(4):047510, 2012.

\bibitem{cameron2014brief}
D.~E. Cameron, C.~J. Bashor, and J.~J. Collins.
\newblock A brief history of synthetic biology.
\newblock {\em Nature Reviews Microbiology}, 12(5):381--390, 2014.

\bibitem{del2017blueprint}
D.~Del~Vecchio, H.~Abdallah, Y.~Qian, and J.~J. Collins.
\newblock A blueprint for a synthetic genetic feedback controller to reprogram
  cell fate.
\newblock {\em Cell systems}, 4(1):109--120, 2017.

\bibitem{del2016control}
D.~Del~Vecchio, A.~J. Dy, and Y.~Qian.
\newblock Control theory meets synthetic biology.
\newblock {\em Journal of The Royal Society Interface}, 13(120):20160380, 2016.

\bibitem{del2015biomolecular}
D.~Del~Vecchio and R.~M. Murray.
\newblock {\em Biomolecular feedback systems}.
\newblock Princeton University Press Princeton, NJ, 2015.

\bibitem{didovyk2015distributed}
A.~Didovyk, O.~I. Kanakov, M.~V. Ivanchenko, J.~Hasty, R.~Huerta, and
  L.~Tsimring.
\newblock Distributed classifier based on genetically engineered bacterial cell
  cultures.
\newblock {\em ACS synthetic biology}, 4(1):72--82, 2015.

\bibitem{friedland2009synthetic}
A.~E. Friedland, T.~K. Lu, X.~Wang, D.~Shi, G.~Church, and J.~J. Collins.
\newblock Synthetic gene networks that count.
\newblock {\em science}, 324(5931):1199--1202, 2009.

\bibitem{gardner2000construction}
T.~S. Gardner, C.~R. Cantor, and J.~J. Collins.
\newblock Construction of a genetic toggle switch in escherichia coli.
\newblock {\em Nature}, 403(6767):339--342, 2000.

\bibitem{hasnain2019optimal}
A.~Hasnain, N.~Boddupalli, and E.~Yeung.
\newblock Optimal reporter placement in sparsely measured genetic networks
  using the koopman operator.
\newblock {\em arXiv preprint arXiv:1906.00944}, 2019.

\bibitem{hasnain2019data}
A.~Hasnain, S.~Sinha, Y.~Dorfan, A.~E. Borujeni, Y.~Park, P.~Maschhoff,
  U.~Saxena, J.~Urrutia, N.~Gaffney, D.~Becker, et~al.
\newblock A data-driven method for quantifying the impact of a genetic circuit
  on its host.
\newblock In {\em 2019 IEEE Biomedical Circuits and Systems Conference
  (BioCAS)}, pages 1--4. IEEE, 2019.

\bibitem{jacob1961genetic}
F.~Jacob and J.~Monod.
\newblock Genetic regulatory mechanisms in the synthesis of proteins.
\newblock {\em Journal of molecular biology}, 3(3):318--356, 1961.

\bibitem{jayanthi2013retroactivity}
S.~Jayanthi, K.~S. Nilgiriwala, and D.~Del~Vecchio.
\newblock Retroactivity controls the temporal dynamics of gene transcription.
\newblock {\em ACS synthetic biology}, 2(8):431--441, 2013.

\bibitem{kaiser2019data}
E.~Kaiser, J.~N. Kutz, and S.~L. Brunton.
\newblock Data-driven approximations of dynamical systems operators for
  control.
\newblock {\em arXiv preprint arXiv:1902.10239}, 2019.

\bibitem{korda2018linear}
M.~Korda and I.~Mezic.
\newblock Linear predictors for nonlinear dynamical systems: Koopman operator
  meets model predictive control.
\newblock {\em Automatica}, 93:149--160, 2018.

\bibitem{korda2019optimal}
M.~Korda and I.~Mezic.
\newblock Optimal construction of koopman eigenfunctions for prediction and
  control.
\newblock {\em hal-02278835}, 2019.

\bibitem{lusch2018deep}
B.~Lusch, J.~N. Kutz, and S.~L. Brunton.
\newblock Deep learning for universal linear embeddings of nonlinear dynamics.
\newblock {\em Nature communications}, 9(1):4950, 2018.

\bibitem{mauroy2019koopman}
A.~Mauroy and J.~Goncalves.
\newblock Koopman-based lifting techniques for nonlinear systems
  identification.
\newblock {\em IEEE Transactions on Automatic Control}, 2019.

\bibitem{mezic2005spectral}
I.~Mezic.
\newblock Spectral properties of dynamical systems, model reduction and
  decompositions.
\newblock {\em Nonlinear Dynamics}, 41(1-3):309--325, 2005.

\bibitem{millacura2019parallel}
F.~A. Millacura, B.~Largey, and C.~E. French.
\newblock Parallel: a novel population-based approach to biological logic
  gates.
\newblock {\em Frontiers in Bioengineering and Biotechnology}, 7:46, 2019.

\bibitem{nielsen2016genetic}
A.~A. Nielsen, B.~S. Der, J.~Shin, P.~Vaidyanathan, V.~Paralanov, E.~A.
  Strychalski, D.~Ross, D.~Densmore, and C.~A. Voigt.
\newblock Genetic circuit design automation.
\newblock {\em Science}, 352(6281):aac7341, 2016.

\bibitem{otto2019linearly}
S.~E. Otto and C.~W. Rowley.
\newblock Linearly recurrent autoencoder networks for learning dynamics.
\newblock {\em SIAM Journal on Applied Dynamical Systems}, 18(1):558--593,
  2019.

\bibitem{paun2005dna}
G.~Paun, G.~Rozenberg, and A.~Salomaa.
\newblock {\em DNA computing: new computing paradigms}.
\newblock Springer Science \& Business Media, 2005.

\bibitem{proctor2016dynamic}
J.~L. Proctor, S.~L. Brunton, and J.~N. Kutz.
\newblock Dynamic mode decomposition with control.
\newblock {\em SIAM Journal on Applied Dynamical Systems}, 15(1):142--161,
  2016.

\bibitem{proctor2018generalizing}
J.~L. Proctor, S.~L. Brunton, and J.~N. Kutz.
\newblock Generalizing koopman theory to allow for inputs and control.
\newblock {\em SIAM Journal on Applied Dynamical Systems}, 17(1):909--930,
  2018.

\bibitem{rowley2009spectral}
C.~W. Rowley, I.~Mezic, S.~Bagheri, P.~Schlatter, and D.~S. Henningson.
\newblock Spectral analysis of nonlinear flows.
\newblock {\em Journal of Fluid Mechanics}, 641:115, 2009.

\bibitem{sarpeshkar2014analog}
R.~Sarpeshkar.
\newblock Analog synthetic biology.
\newblock {\em Philosophical Transactions of the Royal Society A: Mathematical,
  Physical and Engineering Sciences}, 372(2012):20130110, 2014.

\bibitem{schmid2010dynamic}
P.~J. Schmid.
\newblock Dynamic mode decomposition of numerical and experimental data.
\newblock {\em Journal of fluid mechanics}, 656:5--28, 2010.

\bibitem{shah2018reprogramming}
R.~Shah and D.~Del~Vecchio.
\newblock Reprogramming cooperative monotone dynamical systems.
\newblock In {\em Submitted to IEEE Conference on Decision and Control}, 2018.

\bibitem{shin2020programming}
J.~Shin, S.~Zhang, B.~S. Der, A.~A. Nielsen, and C.~A. Voigt.
\newblock Programming escherichia coli to function as a digital display.
\newblock {\em Molecular Systems Biology}, 16(3), 2020.

\bibitem{singh2014recent}
V.~Singh.
\newblock Recent advances and opportunities in synthetic logic gates
  engineering in living cells.
\newblock {\em Systems and synthetic biology}, 8(4):271--282, 2014.

\bibitem{sinha2016operator}
S.~Sinha, U.~Vaidya, and R.~Rajaram.
\newblock Operator theoretic framework for optimal placement of sensors and
  actuators for control of nonequilibrium dynamics.
\newblock {\em Journal of Mathematical Analysis and Applications},
  440(2):750--772, 2016.

\bibitem{sontag2004inferring}
E.~Sontag, A.~Kiyatkin, and B.~N. Kholodenko.
\newblock Inferring dynamic architecture of cellular networks using time series
  of gene expression, protein and metabolite data.
\newblock {\em Bioinformatics}, 20(12):1877--1886, 2004.

\bibitem{tabor2009synthetic}
J.~J. Tabor, H.~M. Salis, Z.~B. Simpson, A.~A. Chevalier, A.~Levskaya, E.~M.
  Marcotte, C.~A. Voigt, and A.~D. Ellington.
\newblock A synthetic genetic edge detection program.
\newblock {\em Cell}, 137(7):1272--1281, 2009.

\bibitem{takeishi2017learning}
N.~Takeishi, Y.~Kawahara, and T.~Yairi.
\newblock Learning koopman invariant subspaces for dynamic mode decomposition.
\newblock In {\em Advances in Neural Information Processing Systems}, pages
  1130--1140, 2017.

\bibitem{ward2009comparison}
C.~Ward, E.~Yeung, T.~Brown, B.~Durtschi, S.~Weyerman, R.~Howes, J.~Goncalves,
  H.~Sandberg, and S.~Warnick.
\newblock A comparison of network reconstruction methods for chemical reaction
  networks.
\newblock In {\em Proceedings of the Foundations for Systems Biology and
  Engineering Conference}, pages 197--200, 2009.

\bibitem{williams2016extending}
M.~O. Williams, M.~S. Hemati, S.~T. Dawson, I.~G. Kevrekidis, and C.~W. Rowley.
\newblock Extending data-driven koopman analysis to actuated systems.
\newblock {\em IFAC-PapersOnLine}, 49(18):704--709, 2016.

\bibitem{yeung2017biophysical}
E.~Yeung, A.~J. Dy, K.~B. Martin, A.~H. Ng, D.~Del~Vecchio, J.~L. Beck, J.~J.
  Collins, and R.~M. Murray.
\newblock Biophysical constraints arising from compositional context in
  synthetic gene networks.
\newblock {\em Cell systems}, 5(1):11--24, 2017.

\bibitem{yeung2015global}
E.~Yeung, J.~Kim, J.~Gon{\c{c}}alves, and R.~M. Murray.
\newblock Global network identification from reconstructed dynamical structure
  subnetworks: Applications to biochemical reaction networks.
\newblock In {\em 2015 54th IEEE Conference on Decision and Control (CDC)},
  pages 881--888. IEEE, 2015.

\bibitem{yeung2019learning}
E.~Yeung, S.~Kundu, and N.~Hodas.
\newblock Learning deep neural network representations for koopman operators of
  nonlinear dynamical systems.
\newblock In {\em 2019 American Control Conference (ACC)}, pages 4832--4839.
  IEEE, 2019.

\bibitem{yeung2018koopman}
E.~Yeung, Z.~Liu, and N.~O. Hodas.
\newblock A koopman operator approach for computing and balancing gramians for
  discrete time nonlinear systems.
\newblock In {\em 2018 Annual American Control Conference (ACC)}, pages
  337--344. IEEE, 2018.

\end{thebibliography}

\end{document}